\documentclass[amstex,12pt,russian,amssymb]{article}

\usepackage{mathtext}
\usepackage[cp1251]{inputenc}
\usepackage[T2A]{fontenc}
\usepackage[russian]{babel}
\usepackage[dvips]{graphicx}
\usepackage{amsmath}
\usepackage{amssymb}
\usepackage{amsxtra}
\usepackage{latexsym}
\usepackage{ifthen}

\begin{document}

\newcounter{lemma}
\newcommand{\lemma}{\par \refstepcounter{lemma}%
{\bf Лема \arabic{lemma}.}}

\newcounter{corollary}
\newcommand{\corollary}{\par \refstepcounter{corollary}%
{\bf Наслідок \arabic{corollary}.}}

\newcounter{remark}
\newcommand{\remark}{\par \refstepcounter{remark}%
{\bf Зауваження \arabic{remark}.}}

\newcounter{theorem}
\newcommand{\theorem}{\par \refstepcounter{theorem}%
{\bf Теорема \arabic{theorem}.}}

\newcounter{proposition}
\newcommand{\proposition}{\par \refstepcounter{proposition}%
{\bf Твердження \arabic{proposition}.}}

\newcounter{example}
\newcommand{\example}{\par \refstepcounter{example}%
{\bf Приклад \arabic{example}.}}

\renewcommand{\refname}{\centerline{\bf Список літератури}}

\renewcommand{\figurename}{Мал.}

\newcommand{\proof}{{\it Доведення.\,\,}}

\noindent УДК 517.5

\medskip\medskip
{\bf О.П.~Довгопятий} (Житомирський державний університет імені
Івана Фран\-ка)

{\bf Є.О.~Севостьянов} (Житомирський державний університет імені
Івана Фран\-ка; Інститут прикладної математики і механіки НАН
України, м.~Слов'янськ)

\medskip\medskip
{\bf O.P.~Dovhopiatyi} (Zhytomyr Ivan Franko State University)

{\bf E.A.~Sevost'yanov} (Zhytomyr Ivan Franko State University;
Institute of Applied Ma\-the\-ma\-tics and Mechanics of NAS of
Ukraine, Slov'yans'k)

\medskip\medskip\medskip
{\bf Про оцінки спотворення відображень з умовою Полецького в
областях з нерівністю Пуанкаре}

{\bf On distortion estimates of mappings with the Poletsky condition
in domains with Poincare inequality}

\medskip\medskip\medskip\medskip
Стаття присвячена дослідженню відображень, які спотворюють модуль
сімей кривих по типу нерівності Полецького. У межових точках області
отримана оцінка спотворення відстані для таких відображень за умови,
що їх характеристика має скінченні інтегральні середні по кулях, а
область значення відображення є регулярною за Альфорсом областю з
нерівністю Пуанкаре. У рукопису розглянуті випадки гомеоморфізмів і
відображень з розгалуженням, області визначення з гарними межами та
області з простими кінцями.

\medskip\medskip
The article is devoted to the study of mappings that distort the
modulus of families of paths according to the Poletsky inequality
type. At the boundary points of the domain, we have obtained an
estimate of the distance distortion for such mappings provided that
their characteristic has finite integral averages over the balls,
and the image domain of mappings is Ahlfors regular and such that
the Poincare inequality holds. The manuscript considers the cases of
homeomorphisms and mappings with branching, as well as definition
domains with good boundaries and domains with prime ends.

\newpage
{\bf 1. Вступ.} Добре відомо, що квазіконформні відображення і
відображення зі скінченним спотворенням за певних умов задовольняють
нерівність Гельдера. Питання виконання цієї нерівності було
досліджено як у внутрішніх, так і межових точках області різними
авторами, див., напр., \cite{MRV$_2$}, \cite{M}, \cite{RS},
\cite{MRSY} і \cite{AM}. Зокрема, у наших роботах проблема виконання
нерівності Гельдера була досліджена у класах Соболєва та
Орліча-Соболєва. Як правило, ми розглядали відображення одиничної
кулі на себе, що з точки зору методології є істотним припущенням, а
сама нерівність доводилася в точках одиничної сфери, див.
\cite{RSS}, \cite{MSS} та \cite{MS}. У даному рукопису ми
продовжуємо дослідження на цю тему. Ми розглядаємо більш загальний
клас відображень, для яких виконується нерівність типу Полецького
відносно $p$-модуля, $n-1<p\leqslant n.$ Областями визначення і
значень є області, які задовольняють певні умови. Зокрема, ви
вимагаємо, щоб область значень була регулярною за Альфорсом, у якій
виконується нерівність Пуанкаре. Доволі широкий клас областей
задовольняє ці вимоги, у тому числі, одинична куля (див., напр.,
\cite[теорема~10.5]{HK}). В свою чергу, область визначення має бути
локально зв'язною на межі <<у певному сенсі>>. Нижче ми розглядаємо
різні варіанти локальної зв'язності, одні з котрих відповідають
спотворенню евклідової відстані (для гарних меж), а інші --
спотворенню у метриці простих кінців (для поганих меж).

\medskip
Тут і далі
$$
B(x_0, r)=\{x\in {\Bbb R}^n: |x-x_0|<r\}\,,\qquad {\Bbb B}^n=B(0,
1)\,,$$
$$S(x_0,r) = \{ x\,\in\,{\Bbb R}^n : |x-x_0|=r\}\,, {\Bbb S}^{n-1}=S(0,
1)\,,$$
$$\Omega_n=m({\Bbb B}^n)\,, \omega_{n-1}={\mathcal H}^{n-1}({\Bbb S}^{n-1})\,,$$
$m$ -- міра Лебега в ${\Bbb R}^n,$ ${\mathcal H}^{n-1}$ --
$(n-1)$-вимірна міра Хаусдорфа,
\begin{equation}\label{eq49***}
A(x_0, r_1,r_2): =\left\{ x\,\in\,{\Bbb R}^n:
r_1<|x-x_0|<r_2\right\}\,,
\end{equation}
а $M_p(\Gamma)$ позначає {\it $p$-модуль сім'ї кривих $\Gamma$}
(див.~\cite{Va}). Для заданих множин $E,$ $F\subset\overline{{\Bbb
R}^n}$ і області $D\subset {\Bbb R}^n$ позначимо через
$\Gamma(E,F,D)$ сім'ю всіх кривих $\gamma:[a,b]\rightarrow
\overline{{\Bbb R}^n}$ таких, що $\gamma(a)\in E,\gamma(b)\in\,F$ і
$\gamma(t)\in D$ при $t \in (a, b).$ Якщо маємо сім'ю $\Gamma$
кривих $\gamma:[a, b]\rightarrow D$ в області $D,$ то під
$f(\Gamma)$ будемо розуміти сім'ю кривих $\{(f\circ\gamma):[a,
b]\rightarrow f(D), \gamma\in \Gamma\}.$ (Якщо $\gamma:[a,
b]\rightarrow \overline{D},$ $\gamma(t)\in D$ при $a<t<b,$ то під
$f(\gamma)$ розуміємо криву $f(\gamma(t)),$ $t\in (a, b).$ Так само
і для відповідних сімей $\Gamma$ кривих $\gamma,$ чиї кінці не
належать $D,$ але належать замиканню $D$).

Нехай $Q:{\Bbb R}^n\rightarrow [0, \infty]$ -- вимірна за Лебегом
функція, рівна нулю зовні $D.$ Розглянемо наступне поняття,
див.~\cite[разд.~7.6]{MRSY}. Будемо говорити, що відображення
$f:D\rightarrow \overline{{\Bbb R}^n}$ є {\it кільцевим
$Q$-відображенням у точці $x_0\in \overline{D}$ відносно
$p$-модуля,} $x_0\ne \infty,$ $p\geqslant 1,$ якщо для деякого
$r_0=r(x_0)>0$ і довільних $0<r_1<r_2<r_0$ виконується умова
\begin{equation}\label{eq3*!!}
 M_p(f(\Gamma(S(x_0, r_1),\,S(x_0, r_2),\,D)))\leqslant
\int\limits_{A} Q(x)\cdot \eta^p(|x-x_0|)\ dm(x)\,, \end{equation}
де $\eta:(r_1,r_2)\rightarrow [0,\infty ]$ може бути обраною
довільною невід'ємною вимірною за Лебегом функцією, яка задовольняє
нерівність
\begin{equation}\label{eq28*}
\int\limits_{r_1}^{r_2}\eta(r)\,dr \geqslant 1\,.
\end{equation}
Згідно \cite[раздел 7.22]{He}, будемо говорити, шо борелева функція
$\rho\colon  D\rightarrow [0, \infty]$ є {\it верхнім градієнтом}
функції $u\colon D\rightarrow {\Bbb R},$ якщо для всіх спрямлюваних
кривых $\gamma,$ які з'єднують точки $x$ і $y\in D,$ виконується
нерівність $|u(x)-u(y)|\leqslant \int\limits_{\gamma}\rho\,|dx|.$
Будемо також говорити, що у $D$ виконується {\it $(1; p)$-нерівність
Пуанкаре,} $p\geqslant 1,$ якщо знайдуться сталі $C\geqslant 1$ і
$\tau>0$ такі, що для кожної кулі $B\subset D,$ довільної обмеженої
неперервної функції $u\colon D\rightarrow {\Bbb R}$ і будь-якого її
верхнього градієнта $\rho$ виконується нерівність
$$\frac{1}{m(B)}\int\limits_{B}|u(x)-u_B|\,dm(x)
\leqslant C\cdot({\rm diam\,}B)\left(\frac{1}{m(\tau B)}
\int\limits_{\tau B}\rho^p(x)\, dm(x)\right)^{1/p}\,,$$
де $u_B:=\frac{1}{m(B)}\int\limits_{B}u(x)\, dm(x).$ Область $D$
будемо називати {\it регулярною за Альфорсом,} якщо для кожного
$x_0\in D,$ деякої сталої $C\geqslant 1$ і довільного $R<{\rm
diam}\,D$
виконується нерівність
$$\frac{1}{C}R^{n}\leqslant m(B(x_0, R)\cap D)\leqslant CR^{n}.$$
Для множин $A, B\subset{\Bbb R}^n$ ми покладаємо, як звично,
$${\rm diam}\,A=\sup\limits_{x, y\in A}|x-y|\,,\quad {\rm dist}\,(A, B)=
\inf\limits_{x\in A,
y\in B}|x-y|\,.$$
Іноді замість ${\rm diam}\,A$ і ${\rm dist}\,(A, B)$ ми також пишемо
$d(A)$ і $d(A, B),$ відповідно.

\medskip
Для заданих чисел $\delta>0$ і $p\geqslant 1,$ фіксованих областей
$D, D^{\,\prime}\subset {\Bbb R}^n,$ $n\geqslant 2,$ $x_0\in
\partial D,$ континуума $A\subset D$ і заданої
функції $Q:D\rightarrow[0, \infty]$ позначимо через $\frak{F}^{p,
x_0}_{Q, A, \delta}(D, D^{\,\prime})$ сім'ю всіх кільцевих
$Q$-гомеоморфізмів $f$ області $D$ на $D^{\,\prime}$ у точці $x_0$
відносно $p$-модуля, які задовольняють умову ${\rm
diam\,}(f(A))\geqslant\delta.$ Справедлива наступна

\medskip
\begin{theorem}\label{th3} {\sl\, Нехай $x_0\in \partial D,$
$x_0\ne\infty,$ $n-1<p\leqslant n,$ $D^{\,\prime}$ є регулярною за
Альфорсом обмеженою областю з $(1; p)$-нерівністю Пуанкаре.
Припустимо, що виконані наступні умови: 1) існує
$r^{\,\prime}_0=r^{\,\prime}_0(x_0)>0$ таке, що множина $B(x_0,
r)\cap D$ є зв'язною при всіх $0<r<r^{\,\prime}_0;$ 2) знайдеться
така стала $0<C(x_0)<\infty,$ що
\begin{equation}\label{eq1D}
\limsup\limits_{\varepsilon\rightarrow
0}\frac{1}{\Omega_n\cdot\varepsilon^n}\int\limits_{B(x_0,
\varepsilon)\cap D}Q(x)\,dm(x)\leqslant C\,.
\end{equation}
Тоді існують $0<\varepsilon_0=\varepsilon_0(x_0, D^{\,\prime}, n, p,
\delta, A)<1$ і $\widetilde{C}=\widetilde{C}(p, n, C,
D^{\,\prime})>0$ такі, що для всіх $x, y\in B(x_0,
\varepsilon_0^2)\cap D$ і всіх $f\in\frak{F}^{p, x_0}_{Q, A,
\delta}(D, D^{\,\prime})$ виконується нерівність
\begin{equation}\label{eq2.4.3}
|f(x)-f(y)|\leqslant \widetilde{C}\cdot\max\left\{
\frac{1}{\log^{n-1}\frac{\varepsilon_0}{|x-x_0|}},
\frac{1}{\log^{n-1}\frac{\varepsilon_0}{|y-x_0|}}\right\}\,.\end{equation}
}
\end{theorem}

\medskip
\begin{corollary}\label{cor1}
{\sl\, В умовах теореми~\ref{th3}, кожне відображення
$f\in\frak{F}^{p, x_0}_{Q, A, \delta}(D, D^{\,\prime})$ має
неперервне продовження в точку $x_0$ і
\begin{equation}\label{eq24B}
|f(x)-f(x_0)|\leqslant
\widetilde{C}\cdot\frac{1}{\log^{n-1}\frac{\varepsilon_0}{|x-x_0|}}\,.
\end{equation}
  }
\end{corollary}

\medskip
\begin{remark}\label{rem1}
Зокрема, умова~1) теореми~\ref{th3} виконується, якщо $D$ -- опукла
область.
\end{remark}

\medskip
У випадку відображень з розгалуженням також виконується деякий
аналог, див. нижче. Нагадаємо, що відображення
$f:D\rightarrow\overline{{\Bbb R}^n}$ області $D\subset{\Bbb R}^n$
на область $D^{\,\prime}\subset\overline{{\Bbb R}^n}$ називається
{\it замкненим}, якщо $C(f, \partial D)\subset
\partial D^{\,\prime},$ де, як звично, $C(f, \partial D)$ --
гранична множина відображення $f$ на $\partial D.$

\medskip
У подальшому,  в розширеному просторі $\overline{{{\Bbb
R}}^n}={{\Bbb R}}^n\cup\{\infty\}$ використовується {\it сферична
(хордальна) метрика} $h(x,y)=|\pi(x)-\pi(y)|,$ де $\pi$ --
стереографічна проекція $\overline{{{\Bbb R}}^n}$  на сферу
$S^n(\frac{1}{2}e_{n+1},\frac{1}{2})$ в ${{\Bbb R}}^{n+1},$ а саме,
$$h(x,\infty)=\frac{1}{\sqrt{1+{|x|}^2}}\,,$$
\begin{equation}\label{eq3C}
\ \ h(x,y)=\frac{|x-y|}{\sqrt{1+{|x|}^2} \sqrt{1+{|y|}^2}}\,, \ \
x\ne \infty\ne y
\end{equation}
(див., напр., \cite[означення~12.1]{Va}).
У подальшому, для множин $A, B\subset \overline{{\Bbb R}^n}$
покладемо
\begin{equation}\label{eq5}
h(A, B)=\inf\limits_{x\in A, y\in B}h(x, y)\,,\quad
h(A)=\sup\limits_{x, y\in A}h(x ,y)\,,
\end{equation}
де $h$ -- хордальная відстань, визначена в~(\ref{eq3C}).

\medskip
Для заданих $\delta>0,$ областей $D, D^{\,\prime}\subset {\Bbb
R}^n,$ $n\geqslant 2,$ точки $x_0\in
\partial D$ і заданої функції $Q:D\rightarrow[0, \infty]$ позначимо
через $\frak{R}^{p, x_0}_{Q, \delta}(D, D^{\,\prime})$ сім'ю всіх
відкритих, дискретних і замкнених кільцевих $Q$-відображень
$f:D\rightarrow D^{\,\prime}$ відносно $p$-модуля у точці $x_0,$
таких, що знайдеться континуум $K_f\subset D^{\,\prime},$ для якого
${\rm diam\,}(K_f)\geqslant \delta$ і $h(f^{\,-1}(K_f),
\partial D)\geqslant \delta>0.$ Справедлива наступна

\medskip
\begin{theorem}\label{th4} {\sl\,
Нехай $x_0\in \partial D,$ $x_0\ne\infty,$ $n-1<p\leqslant n,$
$D^{\,\prime}$ є регулярною за Альфорсом обмеженою областю з $(1;
p)$-нерівністю Пуанкаре. Припустимо, що виконані наступні умови: 1)
існує $r^{\,\prime}_0=r^{\,\prime}_0(x_0)>0$ таке, що множина
$B(x_0, r)\cap D$ є зв'язною при всіх $0<r<r^{\,\prime}_0;$ 2)
знайдеться така стала $0<C=C(x_0)<\infty,$ що
\begin{equation}\label{eq1E}
\limsup\limits_{\varepsilon\rightarrow
0}\frac{1}{\Omega_n\cdot\varepsilon^n}\int\limits_{B(x_0,
\varepsilon)\cap D}Q(x)\,dm(x)\leqslant C\,.
\end{equation}
Тоді існують $0<\varepsilon_0=\varepsilon_0(x_0, D^{\,\prime}, n, p,
\delta, A)<1$ і $\widetilde{C}=\widetilde{C}(p, n, C,
D^{\,\prime})>0$ такі, що для всіх $x, y\in B(x_0,
\varepsilon^2_0)\cap D$ і всіх $f\in\frak{R}^{p, x_0}_{Q, \delta}(D,
D^{\,\prime})$ виконується нерівність
\begin{equation}\label{eq2}
|f(x)-f(y)|\leqslant \widetilde{C}\cdot\max\left\{
\frac{1}{\log^{n-1}\frac{\varepsilon_0}{|x-x_0|}},
\frac{1}{\log^{n-1}\frac{\varepsilon_0}{|y-x_0|}}\right\}\,.\end{equation}
}
\end{theorem}

\medskip
\begin{corollary}\label{cor2}
{\sl\, В умовах теореми~\ref{th4}, відображення $f$ має неперервне
продовження в точку $x_0$ і
$$|f(x)-f(x_0)|\leqslant \widetilde{C}\cdot\frac{1}{\log^{n-1}\frac{\varepsilon_0}{|x-x_0|}}\,.$$
}
\end{corollary}

\medskip
Теореми~\ref{th3} і \ref{th4} мають відповідні аналоги також для
областей з поганими межами, див. нижче. Означення простого кінця,
яке використовується нижче, може бути знайдено в~\cite{KR$_1$}.
Зокрема, будемо говорити, що кінець $K$ є {\it простим}, якщо $K$
містить ланцюг розрізів $\{\sigma_m\}$, таку що
\begin{equation}\label{eqSIMPLE}
\lim\limits_{m\rightarrow\infty}M(\Gamma(C, \sigma_m, D))=0
\end{equation}
для довільного континуума $C$ в $D,$ де $M$ -- модуль сімей кривих
$\Gamma(C, \sigma_m, D).$ Тут і далі $\overline{D}_P$ позначає
поповнення області $D$ її простими кінцями, а
$E_D=\overline{D}_P\setminus D$ -- множина всіх простих кінців в
$D.$ Говоримо, що обмежена область $D$ в ${\Bbb R}^n$ {\it регулярна
(у квазіконформному сенсі)}, якщо $D$ може бути квазіконформно
відображена на область з локально квазіконформною межею, замикання
якої є компактом в ${\Bbb R}^n,$ крім того, кожен простий кінець
$P\subset E_D$ є регулярним. Зауважимо, що замикання
$\overline{D}_P$ регулярної області $D$ є {\it метризовним}, при
цьому, якщо $g:D_0\rightarrow D$ -- квазіконформне відображення
області $D_0$ з локально квазіконформною межею на область $D,$ то
для $x, y\in \overline{D}_P$ покладемо:
\begin{equation}\label{eq1A}
\rho(x, y):=|g^{\,-1}(x)-g^{\,-1}(y)|\,,
\end{equation}
де для $x\in E_D$ елемент $g^{\,-1}(x)$ розуміється як деяка (єдина)
точка межі $D_0,$ коректно визначена з огляду
на~\cite[теорема~4.1]{Na}.

\medskip
{\it Тілом простого кінця $P_0\in E_D$} будемо називати множину
$$I(P_0)=\bigcap\limits_{m=1}\limits^{\infty}\overline{d_m}\,,$$
де $d_m,$ $m=1,2,\ldots ,$ -- деяка спадна послідовність областей
розрізів, що відповідає $P_0.$ Можна показати, що $I(P_0)$ визначено
коректно (іншими словами, $I(P_0)$ не залежить від обраної
послідовності $d_m,$ $m=1,2,\ldots ,$) крім того, $I(P_0)\subset
\partial D$ (див., напр., \cite[Proposition~1]{KR$_1$}.

\medskip
Для заданих $\delta>0$ і $p\geqslant 1,$ областей $D,
D^{\,\prime}\subset {\Bbb R}^n,$ $n\geqslant 2,$ і точки $P_0\in
E_D,$ континуума $A\subset D$ і вимірної за Лебегом функції
$Q:D\rightarrow[0, \infty]$ позначимо через $\frak{F}^{p, P_0}_{Q,
A, \delta}(D, D^{\,\prime})$ сім'ю всіх гомеоморфізмів $f$ області
$D$ на $D^{\,\prime},$ які задовольняють
співвідношення~(\ref{eq3*!!})--(\ref{eq28*}) для кожного $x_0\in
I(P_0)$ (де $I(P_0)$ позначає тіло простого кінця $P_0$) таких, що
${\rm diam\,}(f(A))\geqslant\delta.$ Виконується наступне
твердження.

\medskip
\begin{theorem}\label{th1} {\sl\, Нехай $P_0\in E_D:=\overline{D}_P\setminus D,$ область $D$ є
регулярною (у квазіконформному сенсі), а область $D^{\,\prime}$ є
регулярною за Альфорсом обмеженою областю з $(1; p)$-нерівністю
Пуанкаре, $n-1<p\leqslant n.$ Припустимо, що виконані наступні
умови: 1) для кожного $y_0\in\partial D$ існує
$r^{\,\prime}_0=r^{\,\prime}_0(y_0)>0$ таке, що множина $B(y_0,
r)\cap D$ є скінченно зв'язною при всіх $0<r<r^{\,\prime}_0,$
причому, для кожної компоненти $K$ зв'язності множини $B(y_0, r)\cap
D$ виконана наступна умова: будь-які $x, y\in K$ можна з'єднати
кривою $\gamma:[a, b]\rightarrow {\Bbb R}^n$ такою, що $|\gamma|\in
K\cap \overline{B(y_0, \max\{|x-y_0|, |y-y_0|\})},$
$$|\gamma|=\{y\in {\Bbb R}^n: \exists\,t\in[a, b]: \gamma(t)=y\};$$
2) для кожного $y_0\in I(P_0)$ знайдеться така стала
$0<C=C(y_0)<\infty,$ що
\begin{equation}\label{eq1F}
\limsup\limits_{\varepsilon\rightarrow
0}\frac{1}{\Omega_n\cdot\varepsilon^n}\int\limits_{B(y_0,
\varepsilon)\cap D}Q(x)\,dm(x)\leqslant C\,.
\end{equation}
Тоді для кожного $P\in E_D=\overline{D}_P\setminus D$ існує $y_0\in
\partial D$ таке, що $I(P)=\{y_0\},$ де $I(P)$ позначає тіло простого кінця
$P.$ Крім того, існують $\rho_0=\rho_0(P_0, D^{\,\prime}, n, p,
\delta, A)>0,$ $\varepsilon_0=\varepsilon_0(P_0, D^{\,\prime}, n, p,
\delta, A)>0$ і число $\widetilde{C}=\widetilde{C}(p, n, C,
D^{\,\prime})>0$ такі, що для всіх $f\in\frak{F}^{p, P_0}_{Q, A,
\delta}(D, D^{\,\prime})$ і $x, y\in B_{\rho}(P_0, \rho_0)\cap D,$
$$B_{\rho}(P_0, \rho_0):=\{x\in\overline{D}_P:\rho(x, P_0)<\rho_0\}\,,$$
виконується нерівність
\begin{equation}\label{eq3A}
|f(x)-f(y)|\leqslant \widetilde{C}\cdot\max\left\{
\frac{1}{\log^{n-1}\frac{\varepsilon_0}{|x-x_0|}},
\frac{1}{\log^{n-1}\frac{\varepsilon_0}{|y-x_0|}}\right\}\,,\end{equation}
де $x_0:=I(P_0).$

}
\end{theorem}

\medskip
\begin{corollary}\label{cor3}
{\sl\, В умовах і позначеннях теореми~\ref{th1} відображення $f$ має
неперервне продовження в точку $P_0\in E_D,$ причому
\begin{equation}\label{eq3B} |f(x)-f(P_0)|\leqslant
\frac{\widetilde{C}}{\log^{n-1}\frac{\varepsilon_0}{|x-x_0|}}\,.\end{equation}
 }
\end{corollary}

\medskip
\begin{remark}\label{rem1A}
Зокрема, умова~1) теореми~\ref{th1} виконується, якщо для кожного
$x_0\in\partial D$ існує $r^{\,\prime}_0=r^{\,\prime}_0(x_0)>0$
таке, що множина $B(x_0, r)\cap D$ є скінченно зв'язною при всіх
$0<r<r^{\,\prime}_0,$ причому, для кожна компонента $K$ зв'язності
множини $\overline{B(x_0, r)}\cap D$ є опуклою. Дійсно, нехай $x,
y\in K.$ З'єднаємо $x$ та $y$ відрізком $\gamma$ всередині $K$ (це
можливо, оскільки зв'язна відкрита множина $K$ також є лінійно
зв'язною, див., напр., \cite[наслідок~13.1]{MRSY}). Нехай для
визначеності $|x-x_0|\geqslant |y-x_0|.$ В силу того, що куля
$\overline{B(x, |x-x_0|)}$ є опуклою, то і весь відрізок $\gamma$
належить $\overline{B(x, |x-x_0|)}.$ Тоді $\gamma$ і є бажаною
кривою.
\end{remark}

\medskip
Для заданих чисел $p\geqslant 1$ і $\delta>0,$ областей $D,
D^{\,\prime}\subset {\Bbb R}^n,$ $n\geqslant 2,$ точки $P_0\in E_D,$
і заданої вимірної за Лебегом функції $Q:D\rightarrow[0, \infty]$
позначимо через $\frak{R}^{p, P_0}_{Q, \delta}(D, D^{\,\prime})$
сім'ю всіх відкритих, дискретних і замкнених відображень $f$ області
$D$ на $D^{\,\prime}$ які задовольняють
умови~(\ref{eq3*!!})--(\ref{eq28*}) у кожній точці $x_0\in I(P_0)$
(де $I(P_0)$ позначає тіло простого кінця $P_0$) таких, що
знайдеться континуум $K_f\subset D^{\,\prime}_f,$ для которого ${\rm
diam\,}(K_f)\geqslant \delta$ і $h(f^{\,-1}(K_f),
\partial D )\geqslant \delta>0.$ Виконується наступна теорема.

\medskip
\begin{theorem}\label{th2} {\sl\,
Нехай область $D$ є регулярною, а область $D^{\,\prime}$ є
регулярною за Альфорсом обмеженою областю з $(1; p)$-нерівністю
Пуанкаре, $n-1<p\leqslant n.$ Припустимо, що виконані наступні
умови: 1) для кожного $x_0\in\partial D$ існує
$r^{\,\prime}_0=r^{\,\prime}_0(x_0)>0$ таке, що множина $B(x_0,
r)\cap D$ є скінченно зв'язною при всіх $0<r<r^{\,\prime}_0,$
причому, для кожної компоненти $K$ зв'язності множини
$\overline{B(x_0, r)}\cap D$ виконана наступна умова: будь-які $x,
y\in K$ можна з'єднати кривою $\gamma:[a, b]\rightarrow {\Bbb R}^n$
такою, що $|\gamma|\in K\cap
\overline{B(x_0, \max\{|x-x_0|, |y-x_0|\})},$ %
$$|\gamma|=\{x\in {\Bbb R}^n: \exists\,t\in[a, b]: \gamma(t)=x\};$$
2) для кожного $x_0\in I(P_0)$ знайдеться така стала
$0<C=C(x_0)<\infty,$ що
\begin{equation}\label{eq1G}
\limsup\limits_{\varepsilon\rightarrow
0}\frac{1}{\Omega_n\cdot\varepsilon^n}\int\limits_{B(x_0,
\varepsilon)\cap D}Q(x)\,dm(x)\leqslant C\,,
\end{equation}
де $I(P_0)$ позначає тіло простого кінця $P_0.$
Тоді для кожного $P\in E_D:=\overline{D}_P\setminus D$ існує $y_0\in
\partial D$ таке, що $I(P)=\{y_0\}.$ Крім того, існують $\rho_0=\rho_0(P_0, D^{\,\prime}, n, p,
\delta)>0,$ $\varepsilon_0=\varepsilon_0(P_0, D^{\,\prime}, n, p,
\delta)>0$ і число $\widetilde{C}=\widetilde{C}(p, n, C,
D^{\,\prime})>0$ такі, що для всіх $f\in\frak{R}^{p, P_0}_{Q,
\delta}(D, D^{\,\prime})$ і $x, y\in B_{\rho}(P_0, \rho_0)\cap D$
виконується нерівність
\begin{equation}\label{eq4A}
|f(x)-f(y)|\leqslant \widetilde{C}\cdot\max\left\{
\frac{1}{\log^{n-1}\frac{\varepsilon_0}{|x-x_0|}},
\frac{1}{\log^{n-1}\frac{\varepsilon_0}{|y-x_0|}}\right\}\,,
 \end{equation}
де $x_0:=I(P_0).$
  }
\end{theorem}

\medskip
\begin{corollary}\label{cor4}
{\sl\, В умовах і позначеннях теореми~\ref{th2} відображення $f$ має
неперервне продовження в точку $P_0\in E_D,$ причому
\begin{equation}\label{eq3BD} |f(x)-f(P_0)|\leqslant
\frac{\widetilde{C}}{\log^{n-1}\frac{\varepsilon_0}{|x-x_0|}}\,.
\end{equation}}
\end{corollary}

\medskip
{\bf 2. Допоміжні відомості.} Наступний результат може бути
знайдений в \cite[Pro\-po\-si\-tion~4.7]{AS}.

\medskip
\begin{proposition}\label{pr_2}{\sl\,
Нехай область $D$ є регулярною за Альфорсом, у якій виконується $(1;
p)$-нерівність Пуанкаре, $n-1<p\leqslant n.$ Тоді існує стала $M>0$
така, що для всіх $x\in D,$ $R>0$ і всіх континуумів $E$ і $F$ у
$B(x, R)$ виконується нерівність
$$M_p(\Gamma(E, F, D))\geqslant \frac{1}{M}\cdot\frac{\min\{{\rm diam}\,E, {\rm
diam}\,F\}}{R^{1+p-n}}\,.$$}
\end{proposition}
Нехай $a>0$ і $\varphi\colon [a,\infty)\rightarrow[0,\infty)$~---
неспадна функція, така що при деяких сталих $\gamma>0,$ $T>0$ і всіх
$t\geqslant T$ виконана нерівність
\begin{equation}\label{eq1B}
\varphi(2t)\leqslant \gamma\cdot\varphi(t)\,.
\end{equation}
Будемо називати такі функції {\it функціями, які задовольняють умову
подвоєння}.

Нехай $\varphi\colon [a,\infty)\rightarrow[0,\infty)$ -- функція з
умовою подвоєння, тоді функція
$\widetilde{\varphi}(t):=\varphi(1/t)$ не зростає та визначена на
пів інтервалі  $(0, 1/a].$ Наступне твердження доведено
в~\cite[лема~3.1]{RSS}.

\medskip
\begin{proposition}\label{pr1}
{\sl Нехай $a>0,$ $\varphi\colon [a,\infty)\rightarrow[0,\infty)$ --
неспадна функція з умовою подвоєння~(\ref{eq1B}), $x_0\in {\Bbb
R}^n,$ $n\geqslant 2,$ і нехай $Q:{\Bbb R}^n\rightarrow [0, \infty]$
-- вимірна за Лебегом функція, для якої існує $0<C<\infty$ таке, що
\begin{equation}\label{eq1AA}
\limsup\limits_{\varepsilon\rightarrow
0}\frac{\varphi(1/\varepsilon)}{\Omega_n\cdot\varepsilon^n
}\int\limits_{B(x_0, \varepsilon)}Q(x)\,dm(x)\leqslant C\,.
\end{equation}
Тоді знайдеться $\varepsilon^{\,\prime}_0>0$ таке, що
\begin{equation}\label{eq2B}
\int\limits_{\varepsilon<|x-x_0|<\varepsilon^{\,\prime}_0}
\frac{\varphi(1/|x-x_0|)Q(x)\,dm(x)}{|x-x_0|^n}\leqslant
C_1\cdot\left(\log\frac{1}{\varepsilon}\right)\,,\qquad\varepsilon\rightarrow
0\,,
\end{equation}
де $C_1:=\frac{\gamma C\Omega_n2^n}{\log 2}.$
}
\end{proposition}

\medskip
Нехай $D\subset {\Bbb R}^n,$ $f:D\rightarrow {\Bbb R}^n$ -- відкрите
дискретне відображення, $\beta: [a,\,b)\rightarrow {\Bbb R}^n$ --
крива і нехай $x\in\,f^{\,-1}(\beta(a)).$ Крива $\alpha:
[a,\,c)\rightarrow D$ називається {\it максимальним $f$-підняттям}
кривої $\beta$ з початком у точці $x,$ якщо $(1)\quad
\alpha(a)=x\,;$ $(2)\quad f\circ\alpha=\beta|_{[a,\,c)};$ $(3)$\quad
при кожному $c<c^{\prime}\leqslant b$ не існує такої кривої
$\alpha^{\prime}: [a,\,c^{\prime})\rightarrow D,$ що
$\alpha=\alpha^{\prime}|_{[a,\,c)}$ і $f\circ
\alpha^{\,\prime}=\beta|_{[a,\,c^{\prime})}.$ Виконується наступне
твердження, див.~\cite[лема~3.12]{MRV$_3$}, див.
також~\cite[лема~3.7]{Vu}.

\medskip
\begin{proposition}\label{pr2}
{\sl Нехай $f:D\rightarrow {\Bbb R}^n,$ $n\geqslant 2,$ -- відкрите
дискретне відображення, нехай $x_0\in D,$ і нехай $\beta:
[a,\,b)\rightarrow {\Bbb R}^n$ -- крива, така що $\beta(a)=f(x_0)$
і, або границя $\lim\limits_{t\rightarrow b}\beta(t)$ існує, або
$\beta(t)\rightarrow \partial f(D)$ при $t\rightarrow b.$ Тоді
$\beta$ має максимальне $f$-підняття кривої $\alpha:
[a,\,c)\rightarrow D$ з початком у точці $x_0.$ Якщо
$\alpha(t)\rightarrow x_1\in D$ при $t\rightarrow c,$ то $c=b$ і
$f(x_1)=\lim\limits_{t\rightarrow b}\beta(t).$ В протилежному
випадку, $\alpha(t)\rightarrow \partial D$ при $t\rightarrow c.$}
\end{proposition}

\medskip
{\bf 3. Доведення теореми~\ref{th3}.} Покладемо
$$\varepsilon_0=\min\{\varepsilon^{\,\prime}_0(x_0),
r_0^{\,\prime}, d(x_0, A), 1\}\,,$$
де $\varepsilon^{\,\prime}_0>0$ -- число з твердження~\ref{pr1}, у
якому ми покладаємо $\varphi\equiv 1,$ крім того, $r_0^{\,\prime}$ і
$A$ -- число і континуум з умов теореми, відповідно.
Нехай $x, y\in B(x_0, \varepsilon^2_0)$ і $f\in \frak{F}^{p,
x_0}_{Q, A, \delta}(D, D^{\,\prime}).$ Тоді $x, y\in B(x_0,
\varepsilon_0).$ Без обмеження загальності ми можемо вважати, що
$|x-x_0|\geqslant|y-x_0|.$
За означенням числа $\varepsilon_0,$
\begin{equation}\label{eq18}
A\subset D\setminus B(x_0, \varepsilon_0)\,.
\end{equation}
Оскільки точки $S(x_0, r)\cap D,$ $0<r<r^{\,\prime}_0,$ є досяжними
з області $D$ деякою кривої $\gamma,$ і множина $B(x_0, r)\cap D$
зв'язна при $0<r<r^{\,\prime}_0,$ точки $x$ і $y$ можна з'єднати
кривою $K,$ яка повністю належить до $\overline{B(x_0, |x-x_0|)}$ і
належить $D.$ Нехай $z, w\in f(A)\subset D^{\,\prime}$ і $u, v\in A$
є такими, що
\begin{equation}\label{eq8}
{\rm diam}\,f(A)=|z-w|=|f(u)-f(v)|\geqslant \delta\,.
\end{equation}
Зауважимо, що $f(K)$ і $f(A)$ є континуумами як неперервний образ
континуума при відображенні $f.$ Зафіксуємо довільним чином $y_0\in
D^{\,\prime}.$ Тоді, оскільки область $D^{\,\prime}$ обмежена, існує
$R_0>0$ таке, що $D^{\,\prime}\subset B(y_0, R_0):=B_{R_0}.$ В
такому випадку, за твердженням~\ref{pr_2}
\begin{equation}\label{eq10}
M_p(\Gamma(f(K), f(A), D^{\,\prime}))\geqslant
\frac{1}{M}\cdot\frac{\min\{{\rm diam}\,f(K), {\rm
diam}\,f(A)\}}{R_0^{1+p-n}}\,.
\end{equation}
Зауважимо, що
\begin{equation}\label{eq17}
\Gamma(|K|, A, D)>\Gamma(S(x_0, |x-x_0|), S(x_0, \varepsilon_0),
D)\,.
\end{equation}
Дійсно, нехай $\gamma\in \Gamma(|K|, A, D),$ $\gamma:[0,
1]\rightarrow D,$ $\gamma(0)\in |K|$ і $\gamma(1)\in A.$ Оскільки $
|K|\subset B(x_0, |x-x_0|),$ з огляду на~(\ref{eq18}), ми отримаємо,
що $A\subset D\setminus B(x_0, |x-x_0|).$ Тоді $|\gamma|\cap B (x_0,
|x-x_0|)\ne\varnothing\ne |\gamma|\cap (D\setminus B(x_0,
|x-x_0|)).$ За~\cite[теорема~1.I.5.46]{Ku} існує $t_1\in (0, 1)$
таке, що $\gamma(t_1)\in S(x_0, |x-x_0|).$ Розглянемо криву
$\gamma_1:=\gamma|_{[t_1, 1]}.$ Нагадаємо, що $|K|\subset B(x_0,
\varepsilon_0),$ крім того, завдяки~(\ref{eq18}) $A\subset
D\setminus B(x_0, \varepsilon_0).$ Тоді $\gamma_1\cap B(x_0,
\varepsilon_0)\ne\varnothing\ne |\gamma|\cap (D\setminus B(x_0,
\varepsilon_0)).$ За~\cite[теорема~1.I.5.46]{Ku} існує $t_2\in (0,
t_1)$ таке, що $\gamma_1(t_1)=\gamma(t_1)\in S (x_0,
\varepsilon_0).$ Покладемо $\gamma_2:=\gamma|_{[t_1, t_2]}.$ Тоді
$\gamma_2$ є підкривою $\gamma$ і $\gamma_2\in\Gamma(S(x_0
,|x-x_0|), S(x_0, \varepsilon_0), D).$ Це доводить~(\ref{eq17}). В
цьому випадку, зі міноруванням модуля сімей кривих (див., напр.,
\cite[теорема~6.4]{Va}), з огляду на~(\ref{eq3*!!}) і~(\ref{eq17})
ми отримаємо, що
$$M_p(\Gamma(f(|K|), f(A), D^{\,\prime}))=
M_p(f(\Gamma(|K|, A, D)))\leqslant $$
\begin{equation}\label{eq20}
\leqslant M_p(f(\Gamma(S(x_0, |x-x_0|), S(x_0, \varepsilon_0),
D)))\leqslant\int\limits_{A} Q(x)\cdot \eta^p(|x-x_0|)\, dm(x)\,,
\end{equation}
де $\eta$ -- довільна вимірна за Лебегом функція, що задовольняє
умову~(\ref{eq28*}) при $r_1=|x-x_0|$ і $r_2=\varepsilon_0.$
Покладемо
$$\eta(t):=\begin{cases}\frac{1}
{\left(\log\frac{\varepsilon_0}{|x-x_0|}\right)^{n/p}
t^{n/p}}\,,& t\in(|x-x_0|, \varepsilon_0)\,, \\
0\,,& t\not\in(|x-x_0|, \varepsilon_0)\,.
\end{cases}$$
Зауважимо, що
$$\int\limits_{|x-x_0|}^{\varepsilon_0}\frac{dt}{t\log\frac{\varepsilon_0}{|x-x_0|}}=\frac{1}{\log\frac{\varepsilon_0}{|x-x_0|}}
\cdot\int\limits_{|x-x_0|}^{\varepsilon_0}\frac{dt}{t}=1\,.$$
Тоді за нерівністю Гельдера ми отримаємо, що
$$1=\int\limits_{|x-x_0|}^{\varepsilon_0}\frac{dt}{t\log\frac{\varepsilon_0}{|x-x_0|}}\leqslant
\left(\int\limits_{|x-x_0|}^{\varepsilon_0}\frac{dt}{t^{\frac{n}{p}}\cdot\left(\log\frac{\varepsilon_0}{|x-x_0|}\right)
^{\frac{n}{p}}}\right)^{\frac{p}{n}}\cdot
\left(\varepsilon_0-|x-x_0|\right)^{\frac{n-p}{n}}\leqslant$$
\begin{equation}\label{eq17C}
\leqslant\left(\int\limits_{|x-x_0|}^{\varepsilon_0}\frac{dt}{t^{\frac{n}{p}}\cdot\left(\log\frac{\varepsilon_0}{|x-x_0|}\right)
^{\frac{n}{p}}}\right)^{\frac{p}{n}}=\left(\int\limits_{|x-x_0|}^{\varepsilon_0}\eta(t)\,dt\right)^{\frac{p}{n}}\,.
\end{equation}
З огляду на~(\ref{eq17C}) випливає, що функція $\eta$ задовольняє
умову~(\ref{eq28*}) при $r_1=|x-x_0|$ і $r_2=\varepsilon_0.$ В цьому
випадку, з огляду на~(\ref{eq20}) ми отримаємо, що
\begin{equation}\label{eq21}
M_p(\Gamma(f(|K|), f(A), D^{\,\prime}))\leqslant
\frac{1}{\log^n\frac{\varepsilon_0}{|x-x_0|}}\int\limits_{A(x_0,
|x-x_0|, \varepsilon_0)} \frac{Q(x)}{{|x-x_0|}^n} \,dm(x)\,.
\end{equation}
Оскільки $|x-x_0|<\varepsilon^2_0,$
\begin{equation}\label{eq3}
\log\frac{1}{|x-x_0|}<2\log\frac{\varepsilon_0}{|x-x_0|}\,.
\end{equation}
За обранням $\varepsilon_0,$ з твердженням~\ref{pr1} і з огляду
на~(\ref{eq21}) ми отримаємо, що
\begin{equation}\label{eq22}
M_p(\Gamma(f(|K|), f(A), D^{\,\prime}))\leqslant
\frac{2}{\log^{n-1}\frac{\varepsilon_0}{|x-x_0|}}\frac{C\Omega_n2^n}{\log
2}\,.
\end{equation}
Поєднуючи~(\ref{eq10}) і~(\ref{eq22}), ми отримаємо, що
\begin{equation}\label{eq23}
\frac{1}{M}\cdot\frac{\min\{{\rm diam}\,f(K), {\rm
diam}\,f(A)\}}{R_0^{1+p-n}}\leqslant
\frac{2}{\log^{n-1}\frac{\varepsilon_0}{|x-x_0|}}\frac{C\Omega_n2^n}{\log
2}\,.
\end{equation}
Оскільки $\omega_{n-1}=n\Omega_n,$ останнє співвідношення можна
переписати у вигляді
\begin{equation}\label{eq24}
\frac{1}{M}\cdot\frac{\min\{{\rm diam}\,f(K), {\rm
diam}\,f(A)\}}{R_0^{1+p-n}}\leqslant
\frac{1}{\log^{n-1}\frac{\varepsilon_0}{|x-x_0|}}\frac{C2^{n+1}}{n\log
2}\,.
\end{equation}
З огляду на~(\ref{eq8}), $\min\{{\rm diam}\,f(K), {\rm
diam}\,f(A)\}\geqslant \min\{{\rm diam}\,f(K), \delta\}.$ Тоді
з~(\ref{eq24}) будемо мати:
\begin{equation}\label{eq24BA}
\min\{{\rm diam}\,f(K), \delta\}\leqslant
\frac{MR_0^{1+p-n}}{\log^{n-1}\frac{\varepsilon_0}{|x-x_0|}}\cdot\frac{C2^{n+1}}{n\log
2}\,.
\end{equation}
Зауважимо, що
$\frac{MR_0^{1+p-n}}{\log^{n-1}\frac{\varepsilon_0}{|x-x_0|}}\cdot\frac{C2^{n+1}}{n\log
2}\rightarrow 0$ при $x\rightarrow x_0.$ Тоді існує
$0<\sigma=\sigma(x_0, M, R_0, n, p, \delta)$ таке, що
\begin{equation}\label{eq26}
\frac{MR_0^{1+p-n}}{\log^{n-1}\frac{\varepsilon_0}{|x-x_0|}}\cdot\frac{C2^{n+1}}{n\log
2}<\delta,\qquad \forall\,\, x\in B(x_0, \sigma)\,.
\end{equation}
Нехай $|x-x_0|<\sigma,$ тоді з огляду на~(\ref{eq24BA}) і
(\ref{eq26}), ми отримаємо, що
\begin{equation}\label{eq25B}
|f(x)-f(y)|\leqslant {\rm diam}\,f(K)\leqslant
\frac{1}{\log^{n-1}\frac{\varepsilon_0}{|x-x_0|}}\cdot\frac{MR_0^{1+p-n}C2^{n+1}}{n\log
2}\,.
\end{equation}
Зауважимо, що за нерівністю трикутника $|f(x)-f(y)|\leqslant 2R_0,$
тому при $|x-x_0|\geqslant \sigma$ маємо:
\begin{equation}\label{eq29}
\frac{1}{\log^{n-1}\frac{\varepsilon_0}{|x-x_0|}}\geqslant
\frac{1}{\log^{n-1}\frac{\varepsilon_0}{\sigma}}:=P_0\,,
\end{equation}
отже,
\begin{equation}\label{eq27}
|f(x)-f(y)|\leqslant 2R_0\leqslant \frac{2R_0}{P_0}\cdot
\frac{1}{\log^{n-1}\frac{\varepsilon_0}{|x-x_0|}}\,. \end{equation}
Покладемо
$\widetilde{C}:=\max\left\{\frac{MR_0^{1+p-n}C2^{n+1}}{n\log 2},
 \frac{2R_0}{P_0}\right\}.$ Стала $\widetilde{C}$ залежить тільки від $p,$ $n,$ $C$
і $D^{\,\prime},$ бо $M$ і $R_0$ цілком визначаються по області
$D^{\,\prime}.$ Доведення завершено.~$\Box$

\medskip
{\it Доведення наслідку~\ref{cor1}}. Якщо б відображення $f$ не мало
б границі при $x\rightarrow x_0,$ то ми б побудували не менше двох
послідовностей $x_m\rightarrow x_0$ і $y_m\rightarrow x_0,$
$m\rightarrow\infty,$ таких що $|f(x_m)-f(y_m)|\geqslant \delta>0$
для деякого додатнього $\delta>0$ і всіх $m=1,2,\ldots .$ Але це
суперечить нерівності~(\ref{eq2.4.3}). Отже, границя $f$ при
$x\rightarrow x_0$ існує. Для доведення~(\ref{eq24B}) залишилося
перейти у~(\ref{eq2.4.3}) до границі при $y\rightarrow x_0.$~$\Box$

\medskip
{\bf 4. Доведення теореми~\ref{th4}.} Покладемо
$$\varepsilon_0=\min\{\varepsilon^{\,\prime}_0(x_0),
r_0^{\,\prime}, \delta, 1\}\,,$$
де $\varepsilon^{\,\prime}_0>0$ -- число з твердження~\ref{pr1}, у
якому ми покладемо $\varphi\equiv 1,$ крім того, $r_0^{\,\prime}$ --
число з умов теореми і $\delta$ -- число з означення класу
$\frak{R}^{p, x_0}_{Q, \delta}(D, D^{\,\prime}).$ Покладемо
$\widetilde{\varepsilon_0}:=\varepsilon^2_0.$ Нехай $x, y\in B(x_0,
\widetilde{\varepsilon_0})$ і $f\in \frak{R}^{p, x_0}_{Q, \delta}(D,
D^{\,\prime}).$ Тоді $x, y\in B(x_0, \varepsilon_0).$ Без обмеження
загальності можна вважати, що $|x-x_0|\geqslant|y-x_0|.$ Нехай
$K_{f}\subset D^{\,\prime}$ -- континуум такий, що ${\rm diam\,}
(K_{f})\geqslant \delta$ і $h(f^{\,-1}(K_{f}),
\partial D)\geqslant \delta>0$ (він існує за означенням класу
$\frak{R}^{p, x_0}_{Q, \delta}(D, D^{\,\prime})$). За означенням
$\varepsilon_0,$
\begin{equation}\label{eq18A}
f^{\,-1}(K_{f})\subset D\setminus B(x_0, \varepsilon_0)\,.
\end{equation}
Оскільки точки $S(x_0, r)\cap D,$ $0<r<r^{\,\prime}_0,$ є досяжними
з області $D$ за допомогою деякої кривої $\gamma$ і множина $B(x_0,
r)\cap D$ є зв'язною при $0<r<r^{\,\prime}_0,$ точки $x$ і $y$ можна
з'єднати кривою $K,$ яка повністю належить кулі $\overline{B(x_0,
|x-x_0|)}$ і також належить~$D.$ Нехай $z, w\in K_{f}\subset
D^{\,\prime}$ є такими, що
\begin{equation}\label{eq8A}
{\rm diam}\,K_{f}=|z-w|\geqslant \delta\,.
\end{equation}
Подальше доведення дуже схоже на доведення попередньої
теореми~\ref{th3}. А саме, зауважимо, що $f(K)$ є континуумом як
неперервний образ континуума $K$ при відображенні $f.$ Зафіксуємо
довільним чином $y_0\in D^{\,\prime}.$ Тоді, оскільки область
$D^{\,\prime}$ обмежена, існує $R_0>0$ таке, що $D^{\,\prime}\subset
B(y_0, R_0):=B_{R_0}.$ В такому випадку, за твердженням~\ref{pr_2}
\begin{equation}\label{eq10C}
M_p(\Gamma(f(K), K_{f}, D^{\,\prime}))\geqslant
\frac{1}{M}\cdot\frac{\min\{{\rm diam}\,f(K), {\rm
diam}\,K_{f}\}}{R_0^{1+p-n}}\,.
\end{equation}
Нехай $\Gamma^{\,*}$ -- сім'я $\gamma:[0, 1)\rightarrow D$ усіх
максимальних піднять кривих $\gamma^{\,\prime}:[0, 1]\rightarrow
D^{\,\prime},$ що належать сім'ї $\Gamma=\Gamma(|f(K)|, K_{f},
D^{\,\prime}),$ при відображенні $f$ з початком у $|K|.$ Такі
підняття існують з огляду на твердження~\ref{pr2}. За цим же
твердженням, завдяки замкненості відображення~$f,$ кожна крива
$\gamma\in \Gamma^{\,*}$ має неперервне продовження $\gamma:[0,
1]\rightarrow D$ у точку $b=1.$ Тоді $\gamma(1)\in f^{\,-1}(K_{f}),$
тобто, $\Gamma^{\,*}\subset \Gamma(|K|, f^{\,-1}(K_{f}), D).$

Міркуючи аналогічно доведенню співвідношення~(\ref{eq17}), можна
довести, що
\begin{equation}\label{eq17A}
\Gamma(|K|, f^{\,-1}(K_{f}), D)>\Gamma(S(x_0, |x-x_0|), S(x_0,
\varepsilon_0), D)\,.
\end{equation}
Зауважимо, що $f(\Gamma^{\,*})=\Gamma=\Gamma(|f(K)|, K_{f},
D^{\,\prime}).$ В такому випадку, за принципом мінорування модуля
(див., напр., \cite[теорема~6.4]{Va}), з огляду на~(\ref{eq17A})
і~(\ref{eq3*!!}) ми отримаємо, що
$$M_p(f(\Gamma^{\,*}))=M_p(\Gamma(|f(K)|,
K_{f}, D^{\,\prime}))\leqslant $$
\begin{equation}\label{eq20A}
\leqslant M_p(f(\Gamma(S(x_0, |x-x_0|), S(x_0, \varepsilon_0),
D)))\leqslant\int\limits_{A} Q(x)\cdot \eta^p(|x-x_0|)\ dm(x)\,,
\end{equation}
де $\eta$ -- довільна вимірна за Лебегом функція, яка задовольняє
співвідношення~(\ref{eq28*}) при $r_1=|x-x_0|,$ $r_2=\varepsilon_0.$
Покладемо
$$\eta(t):=\begin{cases}\frac{1}
{\left(\log\frac{\varepsilon_0}{|x-x_0|}\right)^{n/p}
t^{n/p}}\,,& t\in(|x-x_0|, \varepsilon_0)\,, \\
0\,,& t\not\in(|x-x_0|, \varepsilon_0)\,.
\end{cases}$$
За співвідношенням~(\ref{eq17C}) $\eta$ задовольняє
умову~(\ref{eq28*}) при $r_1=|x-x_0|,$ $r_2=\varepsilon_0.$ Тоді
з~(\ref{eq20A}) випливає, що
\begin{equation}\label{eq21A}
M_p(\Gamma(|f(K)|, K_{f}, D^{\,\prime}))\leqslant
\frac{1}{\log^n\frac{\varepsilon_0}{|x-x_0|}}\int\limits_{A(x_0,
|x-x_0|, \varepsilon_0)} \frac{Q(x)}{{|x-x_0|}^p} \,dm(x)\,.
\end{equation}
Оскільки $|x-x_0|<\varepsilon^2_0,$ виконується
співвідношення~(\ref{eq3}). За обранням $\varepsilon_0,$ з
твердженням~\ref{pr1} і з огляду на~(\ref{eq21A}) ми отримаємо, що
\begin{equation}\label{eq22C}
M_p(\Gamma(f(|K|), K_{f}, D^{\,\prime}))\leqslant
\frac{2}{\log^{n-1}\frac{\varepsilon_0}{|x-x_0|}}\frac{C\Omega_n2^n}{\log
2}\,.
\end{equation}
Поєднуючи~(\ref{eq10C}) і~(\ref{eq22C}), ми отримаємо, що
\begin{equation}\label{eq23C}
\frac{1}{M}\cdot\frac{\min\{{\rm diam}\,f(K), {\rm
diam}\,K_{f}\}}{R_0^{1+p-n}}\leqslant
\frac{2}{\log^{n-1}\frac{\varepsilon_0}{|x-x_0|}}\frac{C\Omega_n2^n}{\log
2}\,.
\end{equation}
Оскільки $\omega_{n-1}=n\Omega_n,$ останнє співвідношення можна
переписати у вигляді
\begin{equation}\label{eq24C}
\frac{1}{M}\cdot\frac{\min\{{\rm diam}\,f(K), {\rm
diam}\,K_{f}\}}{R_0^{1+p-n}}\leqslant
\frac{1}{\log^{n-1}\frac{\varepsilon_0}{|x-x_0|}}\frac{C2^{n+1}}{n\log
2}\,.
\end{equation}
З огляду на~(\ref{eq8A}), $\min\{{\rm diam}\,f(K), {\rm
diam}\,K_{f}\}\geqslant \min\{{\rm diam}\,f(K), \delta\}.$ Тоді
з~(\ref{eq24C}) будемо мати:
\begin{equation}\label{eq24D}
\min\{{\rm diam}\,f(K), \delta\}\leqslant
\frac{MR_0^{1+p-n}}{\log^{n-1}\frac{\varepsilon_0}{|x-x_0|}}\cdot\frac{C2^{n+1}}{n\log
2}\,.
\end{equation}
З огляду на~(\ref{eq26}) та~(\ref{eq27}), зі
співвідношення~(\ref{eq24D}) випливає, що
\begin{equation}\label{eq25C}
|f(x)-f(y)|\leqslant {\rm diam}\,f(K)\leqslant \widetilde{C}\cdot
\frac{1}{\log^{n-1}\frac{\varepsilon_0}{|x-x_0|}}\,,
\end{equation}
де $\widetilde{C}:=\max\left\{\frac{MR_0^{1+p-n}C2^{n+1}}{n\log 2},
 \frac{2R_0}{P_0}\right\}.$ Стала
$\widetilde{C}$ залежить тільки від $p,$ $n,$ $C$ і $D^{\,\prime},$
бо $M$ і $R_0$ цілком визначаються по області $D^{\,\prime}.$
Доведення завершено.~$\Box$

\medskip
{\it Доведення} наслідку~\ref{cor2} повністю аналогічно доведенню
наслідку~\ref{cor1}.~$\Box$

\medskip
{\bf 5. Доведення теореми~\ref{th1}.} Для зручності розіб'ємо
доведення на окремі пункти.

\medskip {\bf I.} Нехай $f\in\frak{F}^{p, P_0}_{Q, A, \delta}(D,
D^{\,\prime}).$ Оскільки множина $B(y_0, r)\cap D$ є скінченно
зв'язною при всіх $y_0\in\partial D$ і $0<r<r^{\,\prime}_0(y_0),$
область $D$ є скінченно зв'язною на своїй межі. Отже, область $D$ є
рівномірною (див. \cite[теорема~3.2]{Na$_2$}). Іншими словами, для
кожного $r>0$ існує число $\delta>0$ таке, що нерівність
\begin{equation}\label{eq17***}
M(\Gamma(F^{\,*},F, D))\geqslant \delta
\end{equation}
виконана для всіх континуумів $F, F^*\subset D$ таких, що
$h(F)\geqslant r$ і $h(F^{\,*})\geqslant r.$

\medskip
{\bf II.} Доведемо, що для кожного $P\in E_D$ існує $y_0\in \partial
D$ таке, що $I(P)=\{y_0\}.$ Це твердження доведемо від супротивного,
а саме, припустимо, що існує простий кінець $P\in E_D,$ який містить
дві точки $x,y\in \partial D,$ $x\ne y.$ В такому випадку, існують
принаймні дві послідовності $x_m, y_m\in d_m,$ $m=1,2,\ldots ,$ які
збігаються до $x$ і $y$ при $m\rightarrow\infty,$ відповідно (тут
$d_m$ позначає спадну послідовність областей, утворену деякою
послідовністю розрізів $\sigma_m,$ що відповідає простому кінцю
$P$). З'єднаємо точки $x_m$ і $d_m$ кривою $\gamma_m$ в області
$d_m.$ Оскільки $x\ne y,$ існує $m_0\in {\Bbb N}$ таке, що
$h(\gamma_m)\geqslant d(x,y)/2,$ $m>m_0.$ Оберемо будь-який
невироджений континуум $C\subset D\setminus d_1.$ Тоді з огляду на
рівномірність області~$D$
\begin{equation}\label{eq1}
M(\Gamma(|\gamma_m|, C, D))\geqslant\delta_0>0
\end{equation}
для деякого $\delta_0>0$ і всіх $m>m_0.$ Нерівність~(\ref{eq1})
суперечить означенню простого кінця $P.$ Дійсно, за визначенням
розрізу $\sigma_m$ будемо мати: $\Gamma(|\gamma_m|, C,
D)>\Gamma(\sigma_m, C, D).$ Тоді з огляду на~(\ref{eqSIMPLE}) будемо
мати, що
$$M(\Gamma(|\gamma_m|, C, D))\leqslant M(\Gamma(\sigma_m, C, D))\rightarrow 0$$
при $m\rightarrow\infty.$ Останнє співвідношення
суперечить~(\ref{eq1}). Отже, $I(P)=\{y_0\}$ для деякого $y_0\in
\partial D.$

\medskip
{\bf III.} Залишилося довести співвідношення~(\ref{eq3A}). Нехай
$x_0:=I(P_0)$ і
$$\varepsilon_0=\min\{\varepsilon^{\,\prime}_0(x_0),
r_0^{\,\prime}, {\rm dist}\,(x_0, A), 1\}\,,$$
де $\varepsilon^{\,\prime}_0>0$ -- число з твердження~\ref{pr1},
$\varphi\equiv 1,$ $r_0^{\,\prime}$ -- число з умов теореми.
Оскільки $I(P_0)=\{x_0\},$ знайдеться $\rho_0=\rho_0(P_0,
D^{\,\prime}, n, p, \delta, A)>0$ таке, що $B_{\rho}(P_0,
\rho_0)\subset B(x_0, \widetilde{\varepsilon_0}),$
$\widetilde{\varepsilon_0}=\varepsilon^2_0.$ За означенням
регулярної області, множина $B_{\rho}(P_0, \rho)$ є зв'язною для
достатньо малих $\rho.$ Отже, ми можемо вважати, що $B_{\rho}(P_0,
\rho_0)$ є зв'язним. Нехай $x, y\in B_{\rho}(P_0, \rho_0)$ і
$f\in\frak{F}^{p, P_0}_{Q, A, \delta}(D, D^{\,\prime}).$ Можна
вважати, що $|x-x_0|\geqslant |y-x_0|.$ За означенням
$r_0^{\,\prime}$ точки $x$ і $y$ можна з'єднати кривою $K,$ яка
міститься у кулі $\overline{B(x_0, |x-x_0|)}.$
Нехай $z, w\in f(A)\subset D^{\,\prime}$ і $u, v\in A$ є такими, що
\begin{equation}\label{eq8B}
{\rm diam}\,f(A)=|z-w|=|f(u)-f(v)|\geqslant \delta\,.
\end{equation}
Міркуючи далі аналогічно доведенню теореми~\ref{th3} після
формули~(\ref{eq8}), приходимо до співвідношення
\begin{equation}\label{eq25D}
|f(x)-f(y)|\leqslant {\rm diam}\,f(K)\leqslant
\frac{\widetilde{C}}{\log^{n-1}\frac{\varepsilon_0}{|x-x_0|}}\,,
\end{equation}
де $\widetilde{C}:=\max\left\{\frac{MR_0^{1+p-n}C2^{n+1}}{n\log 2},
\frac{2R_0}{P_0}\right\},$ $M$ -- стала з твердження~\ref{pr_2},
$R_0>0$ -- число таке, що $D^{\,\prime}\subset B(y_0, R_0),$ $y_0$
-- деяка точка області $D^{\,\prime},$ $C>0$ -- стала з
твердження~\ref{pr1}, а $P_0$ визначається по~(\ref{eq29}). Це і
завершує доведення.~$\Box$

\medskip
{\it Доведення наслідку~\ref{cor3}}. За теоремою~\ref{th1}
$I(P_0)=\{x_0\}.$ Тоді при $x\rightarrow P_0$ і $y\rightarrow P_0$
виконано також $x\rightarrow x_0$ і $y\rightarrow x_0.$ Якщо б
відображення $f$ не мало б границі при $x\rightarrow P_0,$ то ми б
побудували не менше двох послідовностей $x_m\rightarrow P_0$ і
$y_m\rightarrow P_0,$ $m\rightarrow\infty,$ таких що
$|f(x_m)-f(y_m)|\geqslant \delta>0$ для деякого додатнього
$\delta>0$ і всіх $m=1,2,\ldots .$ Але це суперечить
нерівності~(\ref{eq3A}). Отже, границя $f$ при $x\rightarrow P_0$
існує. Для доведення~(\ref{eq24B}) залишилося перейти
у~(\ref{eq2.4.3}) до границі при $y\rightarrow P_0.$~$\Box$

\medskip
{\bf 6. Доведення теореми~\ref{th2}}. У зв'язку з тим, що доведення
цієї теореми дуже схоже на усі попередні, обмежимося лише схемою
доведення. Те, що для кожного $P\in E_D$ існує $y_0\in \partial D$
таке, що $I(P)=\{y_0\},$ може бути встановлено так само, як і при
доведенні теореми~\ref{th1}.

Встановимо співвідношення~(\ref{eq4A}). За доведеним вище,
знайдеться $x_0\in
\partial D$ таке, що $I(P_0)=\{x_0\}.$ Нехай
$$\varepsilon_0=\min\{\varepsilon^{\,\prime}_0(x_0),
r_0^{\,\prime}, 1\}\,,$$
де $\varepsilon^{\,\prime}_0>0$ -- число з твердження~\ref{pr1},
$\varphi\equiv 1,$ $r_0^{\,\prime}$ -- число з умов теореми. Звідси
випливає, що знайдеться $\rho_0=\rho_0(P_0, D^{\,\prime}, n, p,
\delta)>0$ таке, що $B_{\rho}(P_0, \rho_0)\subset B(x_0,
\widetilde{\varepsilon_0}),$
$\widetilde{\varepsilon_0}=\varepsilon^2_0.$ За означенням
регулярної області, множина $B_{\rho}(P_0, \rho)$ є зв'язною для
достатньо малих $\rho.$ Отже, ми можемо вважати, що $B_{\rho}(P_0,
\rho_0)$ є зв'язним. Нехай $x, y\in B_{\rho}(P_0, \rho_0)$ і
$f\in\frak{R}^{p, P_0}_{Q, \delta}(D, D^{\,\prime}).$ Можна вважати,
що $|x-x_0|\geqslant |y-x_0|.$ За означенням $r_0^{\,\prime}$ точки
$x$ і $y$ можна з'єднати кривою $K,$ яка міститься у кулі
$\overline{B(x_0, |x-x_0|)}.$

Нехай $K_{f}\subset D^{\,\prime}$ -- континуум такий, що ${\rm
diam\,} (K_{f})\geqslant \delta$ і $h(f^{\,-1}(K_{f}),
\partial D)\geqslant \delta>0$ (він існує за означенням класу
$\frak{R}^{p, P_0}_{Q, \delta}(D, D^{\,\prime})$).
Нехай також $z, w\in K_{f}\subset D^{\,\prime}$ є такими, що
\begin{equation}\label{eq8C}
{\rm diam}\,K_{f}=|z-w|\geqslant \delta\,.
\end{equation}
Проводячи міркування аналогічні доведенню теореми~\ref{th4} після
формули~(\ref{eq8A}), приходимо до співвідношення~(\ref{eq25C}), що
і завершує доведення.~$\Box$

\medskip
{\it Доведення наслідку~\ref{cor4}} цілком аналогічно доведенню
наслідку~\ref{cor3}.~$\Box$

КОНТАКТНАЯ ИНФОРМАЦИЯ

\medskip

\noindent{{\bf Олександр Петрович Довгопятий} \\
Житомирський державний університет ім.\ І.~Франко\\
вул. Велика Бердичівська, 40 \\
м.~Житомир, Україна, 10 008 \\
e-mail: alexdov1111111@gmail.com}

\medskip
\noindent{{\bf Євген Олександрович Севостьянов} \\
{\bf 1.} Житомирський державний університет ім.\ І.~Франко\\
кафедра математичного аналізу, вул. Велика Бердичівська, 40 \\
м.~Житомир, Україна, 10 008 \\
{\bf 2.} Інститут прикладної математики і механіки
НАН України, \\
вул.~генерала Батюка, 19 \\
м.~Слов'янськ, Україна, 84 100\\
e-mail: esevostyanov2009@gmail.com}

\end{document}